\documentclass[conference]{IEEEtran}
\IEEEoverridecommandlockouts
\usepackage{cite}
\usepackage{amsmath,amssymb,amsfonts}
\usepackage{algorithmic}
\usepackage{graphicx}
\usepackage{textcomp}
\usepackage{xcolor}
\usepackage{subfig}
\usepackage[utf8]{inputenc}

\def\BibTeX{{\rm B\kern-.05em{\sc i\kern-.025em b}\kern-.08em
    T\kern-.1667em\lower.7ex\hbox{E}\kern-.125emX}}

\newcommand{\mj}{\mathrm{j}} 
\newcommand{\dif}{\mathrm{d}}
\newcommand{\mat}[1]{\mathrm{\mathbf{{#1}}}} 
\renewcommand{\vec}[1]{\mbox{\boldmath$#1$}} 
 
\newcommand{\vecphas}[1]{\mbox{\underline{\boldmath$#1$}}}

\newcommand{\grad}{\mathrm{grad}\ }
\newcommand{\curl}{\mathrm{curl}\ }

\begin{document}
\title{Artap: Robust Design Optimization Framework for Engineering Applications}

\author{\IEEEauthorblockN{David P\'anek, Tam\'as Orosz, Pavel Karban}
\IEEEauthorblockA{\textit{Department of theory of electrical engineering} \\
\textit{University of West Bohemia}\\
Pilsen, Czech Republic \\
\{panek50, tamas, karban\}@kte.zcu.cz}
}

\maketitle

\begin{abstract}

The main goal of the Artap project is to provide an extensive infrastructure for robust design optimization, where usually many different numerical solvers have to be used together and the impact of the manufacturing uncertainties have to be minimized. Artap is an open-source software platform, developed jointly with the coupled numerical field solver, Agros-suite. Artap ensures interfaces for a broad collection of optimization algorithms (genetic and evolutionary algorithms, various interfaces to libraries such as Nlopt, Bayesopt, etc .),  tools for machine learning (neural networks, Gaussian processes, etc. ), finite element solvers (Agros-suite, Comsol, Multiphysics, Deal.II). The implemented tools offers an easy and straightforward solution not only for robust design optimization but parameter identification, model order reduction, and shape optimization, as well. Moreover, Artap provides automatic parallelization of the optimization process. The paper presents the structure of the framework and technologies powering the project. The main features of Artap are demonstrated on an induction brazing process design tasks.

\end{abstract}
\vspace{3mm}
\let\thefootnote\relax\footnotetext{Presented in ICDS2019, the paper will be published in IEEE Xplore Database}

\section{Introduction}

Design of electric devices and electrical machines is a complex optimization task, where several physical fields and economical factors have to be considered simultaneously. 
Moreover, real world optimization problems are generally large scale problems with many uncertainties \cite{Ben-Tal2002, li2005multi}. 
Those designs are preferred in the practice, which can be realized relatively easily and not sensitive for the small changes of the design parameters. 
There are several systematic design methodology exist in the literature, which goal is to ensure the validity, replicability and robustness of the final product \cite{Ben-Tal2002, tsui1992overview, rdo}. 
The development of Artap is motivated by a similar, multidisciplinary design optimization task: an induction brazing process development.
The goal of the project was to minimize the number of the waste products and the process time during a mass manufacturing process . 
Because of the complexity and the multidisciplinary manner of the optimized process, many sub-design, parameter identification, model order reduction, control and finite element modelling problems should be solved together to reach a robust solution \cite{ih_review, induction_brazing_1,induction_brazing_2}. 
The accurate solution of these type of problems requires a framework, which can handle the different numerical solvers together. 
Moreover, the key to solve robust optimization problems is to find the most suitable optimization strategy for the problem or the modeled sub-problems \cite{moulard2014software, Ben-Tal2002, rdo}. 
It is hard to chose the right optimization strategy and formulate the problem in an appropriate manner in the beginning of the design task \cite{moulard2014software, moulard2013roboptim}.
Instead of a specific optimization solver, a more general solver is advantageous, because it allows to change the design strategy during the design process and provides a high level abstraction layer with the sufficient interfaces to facilitate the implementation of robust design optimization problems \cite{moulard2014software, biscani2010global}.

There are many similar, open-source optimization framework published \cite{biscani2010global, moulard2013roboptim, moulard2014software, Andersson2019, noauthor_openmdao.org_nodate, noauthor_platypus_nodate, johnson2014nlopt, jones2014scipy}, main concepts of Artap is mainly inspired by Platypus \cite{noauthor_platypus_nodate} and OpenMDAO \cite{noauthor_openmdao.org_nodate}. 
Platypus is a high quality code, written in Python, for evolutionary based optimizations, while OpenMDAO is a high-performance computing platform for systems analysis and multidisciplinary optimization. 
However, none of them offers aid in linking together different software parts and interfaces for both of gradient-free and evolutionary based optimization methods, meta-modeling techniques, surrogate models and FEM solvers as Artap for the purpose of combined analyses \cite{biscani2010global, moulard2013roboptim, moulard2014software, Andersson2019, noauthor_openmdao.org_nodate, noauthor_platypus_nodate, johnson2014nlopt, jones2014scipy, toth2018dual, bangerth2007deal}.

\section{Artap}

\begin{figure}[!tb]
	\centering
	\includegraphics[width= 0.9\linewidth]{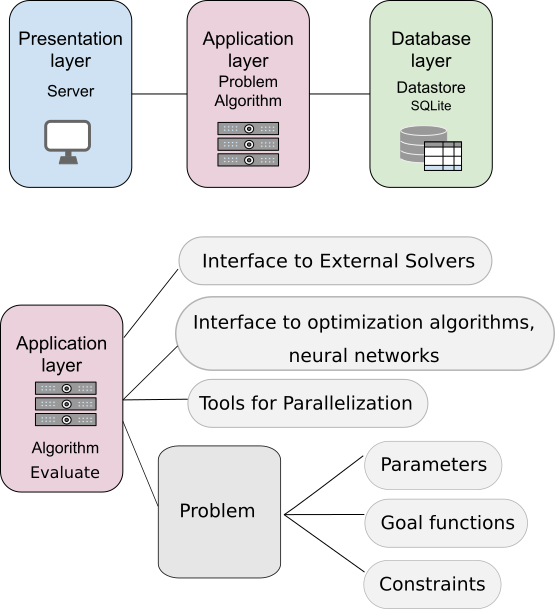}
	\caption{Schematic view of the multi-layered architecture of Artap, which ensures to change the applied optimization algorithm, the applied interfaces without any changes in the problem definition class.}
	\label{fig:aratp_structure}
\end{figure}

Artap is an optimization framework for robust design optimization. 
It is developed within the Department of Theoretical Electrical Engineering in University of West Bohemia jointly with Agros-suite \cite{karban_numerical_2013}. 
Written in Python, Artap aims to facilitate the solution of real-life engineering design problems. 
Where, generally many different numerical solvers have to be used together and the impact of the uncertainties has to be considered during the design optimization process \cite{Ben-Tal2002, li2005multi}. 
A multi-layered architecture is designed to the code, where the problem definition and the optimization solvers, interfaces to other numerical codes are handled separately (Fig.\ref{fig:aratp_structure}). 
Because, the selection of the appropriate optimization algorithm is the key to solve an optimization task in timely manner. However, it is hard task to find the most appropriate solver in the beginning of the task. 

The core of the framework is represented by the application layer which consists of two main classes \texttt{Problem} and \texttt{Algorithm}. The \texttt{Problem} class contains the description of the optimization problem: the constraints, optimized parameters and the goal functions (Fig. \ref{fig:aratp_structure}). 
The \texttt{Algorithm} class contains the selected optimization solver, surrogate model, paralellization settings or links to other external solvers. 
In the simplest case, the user has to redefine only the \texttt{evaluate} function of the \texttt{Algorithm} class, with the actual optimization problem. 
The application layer can be connected to a database through a simple interface and the application layer can be accessed from a web based interface.
A simplified overview of the realized layout is shown in Fig.\ref{fig:aratp_structure}.

At the time of this writing, Artap offers the following optimization algorithms and interfaces:

\begin{itemize}

\item{global and local optimization algorithms coded directly  within  Artap,  including evolutionary and genetic algorithms, like  NSGA-II \cite{deb2002fast, coello2007evolutionary},  $\epsilon$ - MOEA \cite{mishra2002fast},  particle swarm  optimization \cite{kennedy1995eberhart}.}

\item{wrappers for algorithms from the SciPy library \cite{jones2014scipy},  including algorithms for unconstrained optimization: Nelder-Mead method \cite{gao2012implementing}, Powell method \cite{powell1964efficient}, nonlinear conjugate gradient algorithm \cite{nocedal2006sequential},  dog-leg trust-region algorithm \cite{nocedal2006sequential}, BFGS \cite{nocedal2006sequential}, Newton-CG \cite{nocedal2006sequential}, Newton GLTR trust-region algorithm \cite{lenders2018trlib}, trust-region method for unconstrained minimization \cite{conn2000trust};}

\item{wrappers  for  algorithms  from  the  NLopt  library \cite{johnson2014nlopt}<}

\item{wrappers  for  Bayesian optimization from BayesOpt library \cite{martinez2014bayesopt}, efficient implementation of the Bayesian optimization methodology for nonlinear optimization, experimental design and hyper-parameter tuning;}

\item{wrappers for Scikit-learn \cite{pedregosa2011scikit}, machine learning algorithms for medium-scale supervised and unsupervised problems.}

\end{itemize}

The applicability of Artap framework is tested on several real-life engineering design examples. Three of them is shown in the next section to illustrate the applicability of the code: a parameter identification and model order reduction problem for reduced order modelling of induction brazing, a shape optimization of a rotating electrical machines and solution of an analytically validated, TEAM benchmark problem \cite{di2018benchmark} with a neural network based surrogate model and evolutionary algorithms.

\section{Illustrative Examples}
\subsection{Control of Brazing Process}
\begin{figure}[!ht]
	\centering
	\includegraphics[width=.9 \linewidth]{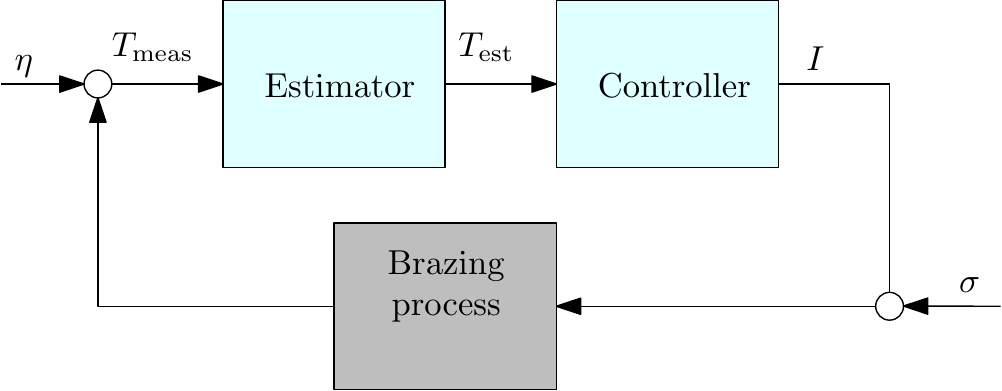}
	\caption{Block diagram of the brazing control method.}
	\label{fig:process_chart}
\end{figure}

\begin{figure}[ht]
	\subfloat[]{\includegraphics[width=.48 \linewidth]{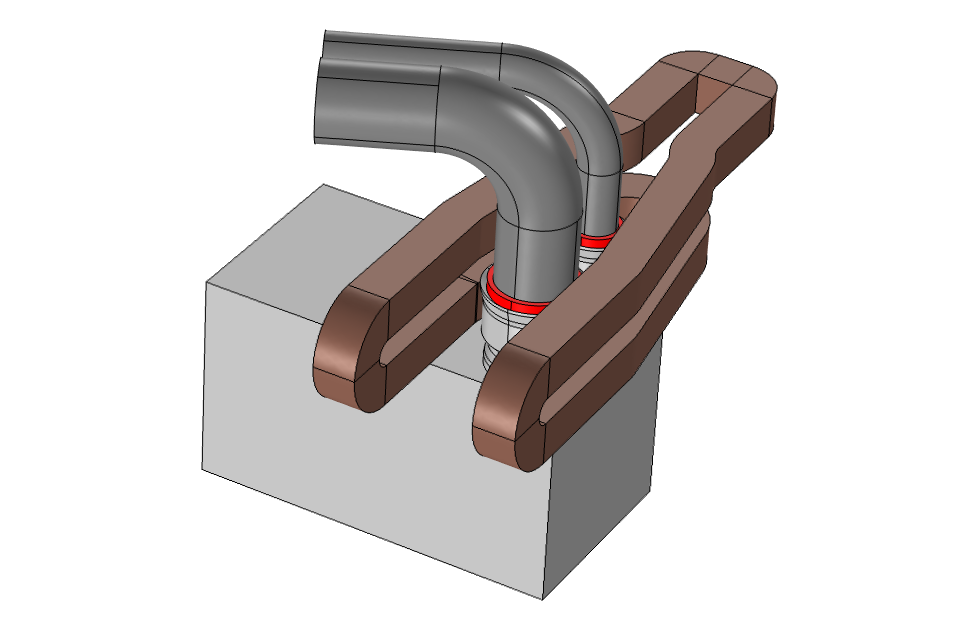} \label{fig:geometry}}
	\subfloat[]{\includegraphics[width=.48 \linewidth]{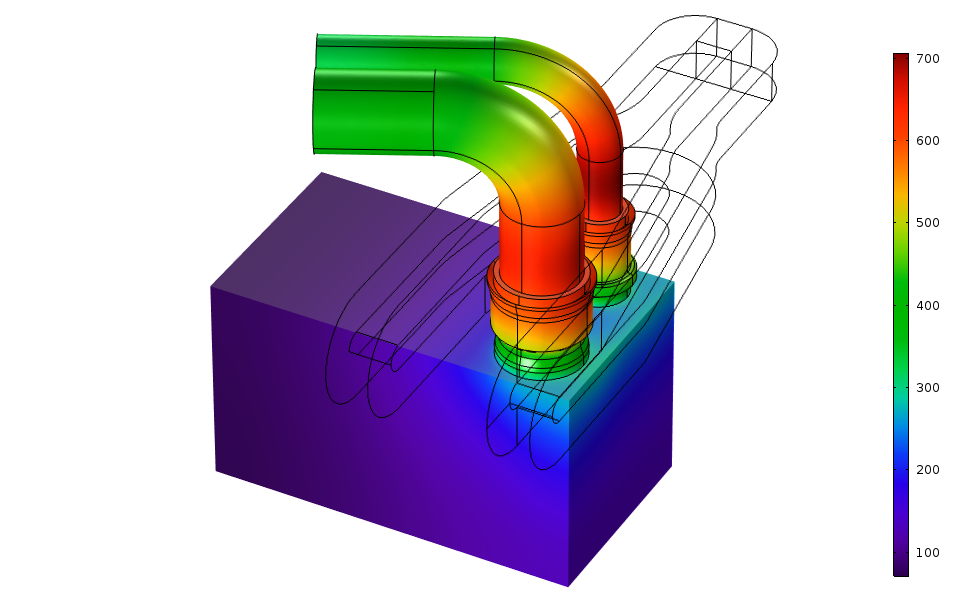} \label{fig:field}}
	\caption{(a) Geometry of the model problem. (b) Distribution of the temperature for $t = 27$[s].}
\end{figure}

Induction brazing of aluminum parts is a relatively fast and clean assembly process. 
However, it is sensitive to the temperature differences and the manufacturing uncertainties. It is necessary to control the process to ensure the required quality \cite{induction_brazing_1, induction_brazing_2}. The 3D model of the inductor and the workpiece (two heat exchanger pipes with sleeves and the evaporator) is shown in Fig.\,\ref{fig:geometry}. The main problem with the control is that the full surface of the brazed components cannot be measured during the brazing process. The following solution is based on a reduced order model of the induction brazing process which can determine any temperature value from the accessible measurement spots. Therefore, a reduced model based estimator is constructed to model the process.
The schematic of the control process is depicted in Fig. \ref{fig:process_chart}. 

The process of induction brazing was modelled by FEM. The geometry of the problem is depicted in Fig. \ref{fig:geometry}. The physical background has been described by the weakly coupled magneto-thermal system, where the magnetic field is described by the Helmholtz-equation:
\begin{equation}
\Delta\vecphas{A}-\mj\cdot\omega\gamma\mu\vecphas{A}=-\mu\vecphas{J}_{\mathrm{e}}\,,
\label{eq:Helmholtz_equation}
\end{equation}
where $\vec{A}$ is the magnetic vector potential, symbol $\omega$ means for frequency, $\gamma$ stands for conductivity, $\mu$ for magnetic permeability and $J_\mathrm{e}$ is the current density which represents the control of the system.
The heat transfer equation is given in the following form: 
\begin{equation}
\nabla \cdot \left(\lambda\nabla T\right)=\rho c_{\rm{p}}\frac{\partial T}{\partial t} - w_{\rm{J}}\,.
\label{eq:heat_transfer_equation}
\end{equation}
The numerical solution of this equation system -- without any simplification -- on an appropriate mesh leads to millions of equations, which cannot be used for real-time control.
There are many possibilities to create a reduced order model with Artap, two possibilities is shown in this paper \cite{induction_brazing_1, induction_brazing_2}. 
Firstly, the problem is formulated as a parameter identification task. In this case the structure of the system is given and parameters of the system are the subject of the optimization. Let the reduced model is described by discrete-time LTI system in the following form:
\begin{eqnarray}
    \vec{x}_{k+1} &=& \vec{A} \vec{x}_k + \vec{B}u_k, \\
     y_{k} &=& \vec{C} \vec{x}_k,
\end{eqnarray}
where $\vec{A}$ matrix represents the discrete-time in Schwartz-form, which can be described in our case, a fourth order system in the following way:
\begin{equation}
  \vec{A} = \left[ \begin{array}{rrrr}
     \Delta_1  &  \delta_1 \Delta_2 &  \delta_1 \delta_2 \Delta_3 &  \delta_1 \delta_2 \delta_3 \Delta_4  \\
     \delta_{1} &   -\Delta_1 \Delta_2 & - \Delta_{1} \delta_{2} \Delta_{3} & - \Delta_{1} \delta_{2} \delta_{3} \Delta_4  \\
     0    & \delta_{2}   &  -\Delta_2 \Delta_3  &  -\Delta_2 \delta_3 \Delta_4  \\
     0    &  0   &  -\delta_3  &  -\Delta_3 \Delta_4  
  \end{array} \right],
\end{equation}
\begin{displaymath}
    \vec{B} = [1, -\Delta_1, -\Delta_2,  -\Delta_3]^\mathrm{T}, \quad \vec{C} = [\gamma_1, \gamma_2, \gamma_3, \gamma_4], 
\end{displaymath}
where $\delta = \sqrt{1 - \Delta^2}$.  
Therefore to describe an n-th order system $2n$ unknown coefficients has to be determined.  
The objective function can be formulated by the following expression:
\begin{equation}
    F(\vec{b}) = \max(\vec{T}_\mathrm{ref} - \vec{T}),
\end{equation}
where $\vec{b} = [\Delta_1, \Delta_2, \Delta_3, \Delta_4, \gamma_1, \gamma_2, \gamma_3, \gamma_4]$ is a vector of the unknown parameters, $T_\mathrm{ref}$ is the temperature vector,it is obtained by the solution of the full field model. $\vec{T}$ is obtained by the usage of the particular variant of reduced system. 

During the second approach, the reduced model is constructed by Proper Orthogonal Decomposition (POD) \cite{chatterjee2000introduction}, which is based on the transformation of the equation system of the full field model. Moreover, it is assumed, that the $U(\vec{x}, t)$ matrix contains all of the FEM solutions. Where the rows represent the time evaluation of a given solution and columns represent the full solution at a selected time.
After the decomposition -- with principal component analysis -- the correlation matrix is $\vec{C}_\mathrm{u}$ is decomposed to 
\begin{equation}
\bar{\vec{C}}_\mathrm{y} = \vec{E} \vec{D} \vec{E}^T.
\end{equation}  
To obtain a reduced model, firstly, the matrices $E$ and $D$ are sorted, according to the size of the eigenvalues. After it, appropriate number of eigenvectors are selected (columns of matrix \vec{E}). Then the reduced model is approximated by the solution of the following expression:
\begin{equation}
\vec{y}_i = \vec{E} \hat{\vec{y}}_i.
\end{equation}
The system of equations for the full field model can be written in the following form:
\begin{equation}
\vec{M} \dot{\vec{y}}(t) + \vec{S} \vec{y}(t) = \vec{F}(t),
\end{equation} 
and the reduced system is then
\begin{equation}
\vec{E}^\mathrm{T}\vec{M} \vec{E} \dot{\hat{\vec{y}}}(t) + \vec{E}^\mathrm{T}\vec{S}\vec{E} \hat{\vec{y}}(t) = \vec{E}^\mathrm{T} \vec{F}(t).
\end{equation} 

Then the system is transformed to the required form by theory of control systems:
\begin{eqnarray}
\vec{\dot{x}}(t) &=& \mathbf{A} \vec{x}(t) + \mat{B} \vec{u}(t)\,,  \\
\vec{y}(t) &=& \mat{C} \vec{x}(t)\,, \nonumber
\label{eq:linear_system}
\end{eqnarray}
where 
$\mat{A} = (\mat{E}^{\rm T} \mat{M} \mat{E})^{-1} \mat{E}^{\rm T} \mat{S} \mat{E}$ is the dynamical matrix of the system, $\mat{B}\vec{u}(t) = (\mat{E}^{\rm T} \mat{M} \mat{E})^{-1} \mat{E}^{\rm T} \mat{F}(t)$ is the input multiplied by the input matrix and $\mat{C}$ is the output matrix representing the way, how the state vector $\vec{x}$ is transformed to a measurable output $\vec{y}(t)$. 

\begin{figure}[!t] 
	\centering
	\includegraphics[width=0.95\linewidth]{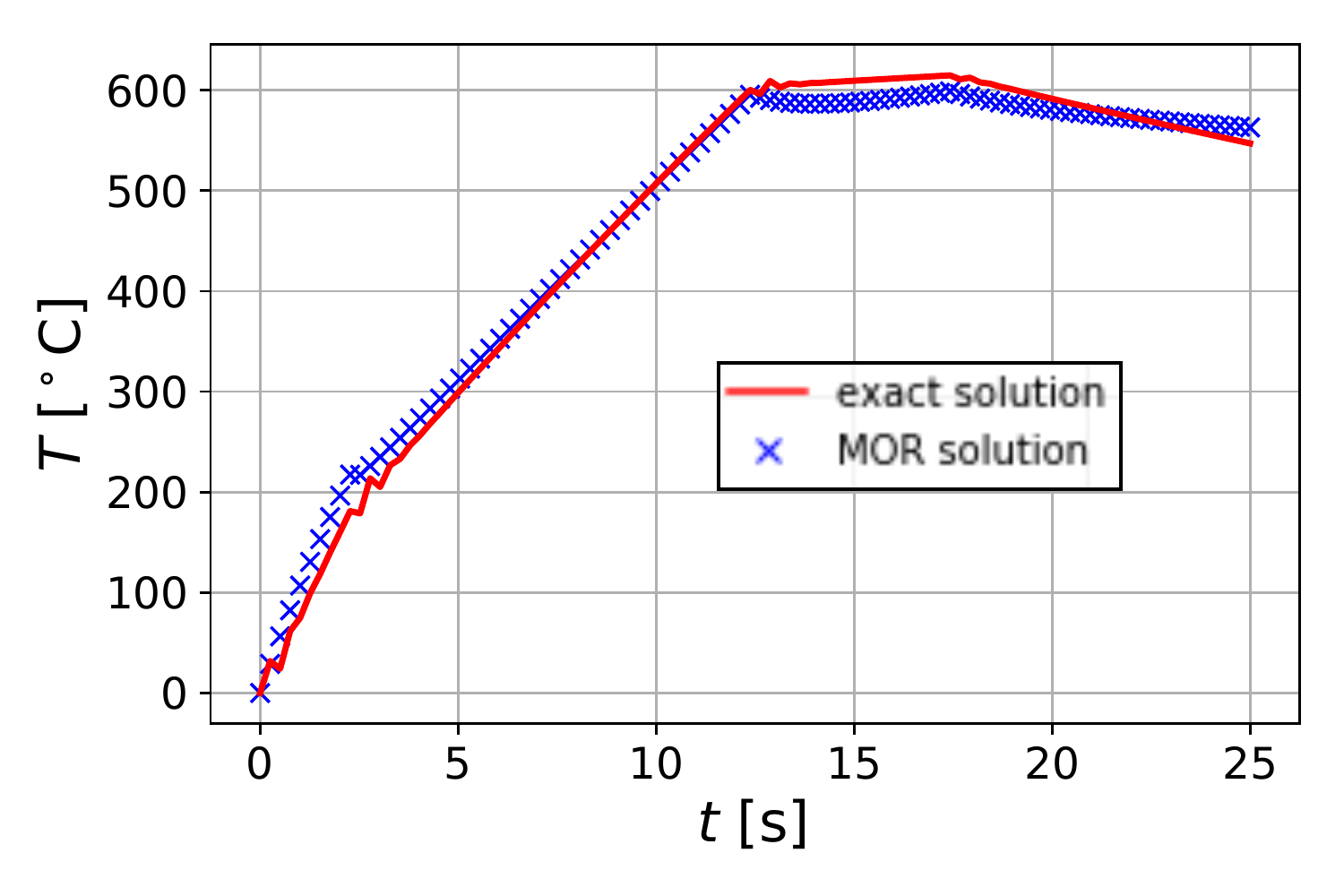}
	\caption{The picture shows the comparison of the results, the time response of the two reduced systems and  the solution of the FEM based full field model. The red line represents the exact solution (Full model), the blue line represents the output of the reduced order model, which was generated by parameters identification and the blue crosses represents the results of the model order reduction.
	}
	\label{fig:identification_mor}
\end{figure}

The resulted dynamic behaviour of the reduced systems is depicted in Fig. \ref{fig:identification_mor}. The results are compared with the full field solution. The results shows that the reduced order models can accurately estimate the temperature value in the non-measurable spots. The difference between solution of the full model and the reduced order model based estimations are lesser than 1\%. Due to the simplicity and the accuracy of these solutions, they can be used for on-line control of induction brazing processes.

\subsection{Reluctance Motor Optimization}
\begin{figure}[!ht]
	\centering
	\includegraphics[width=.7 \linewidth]{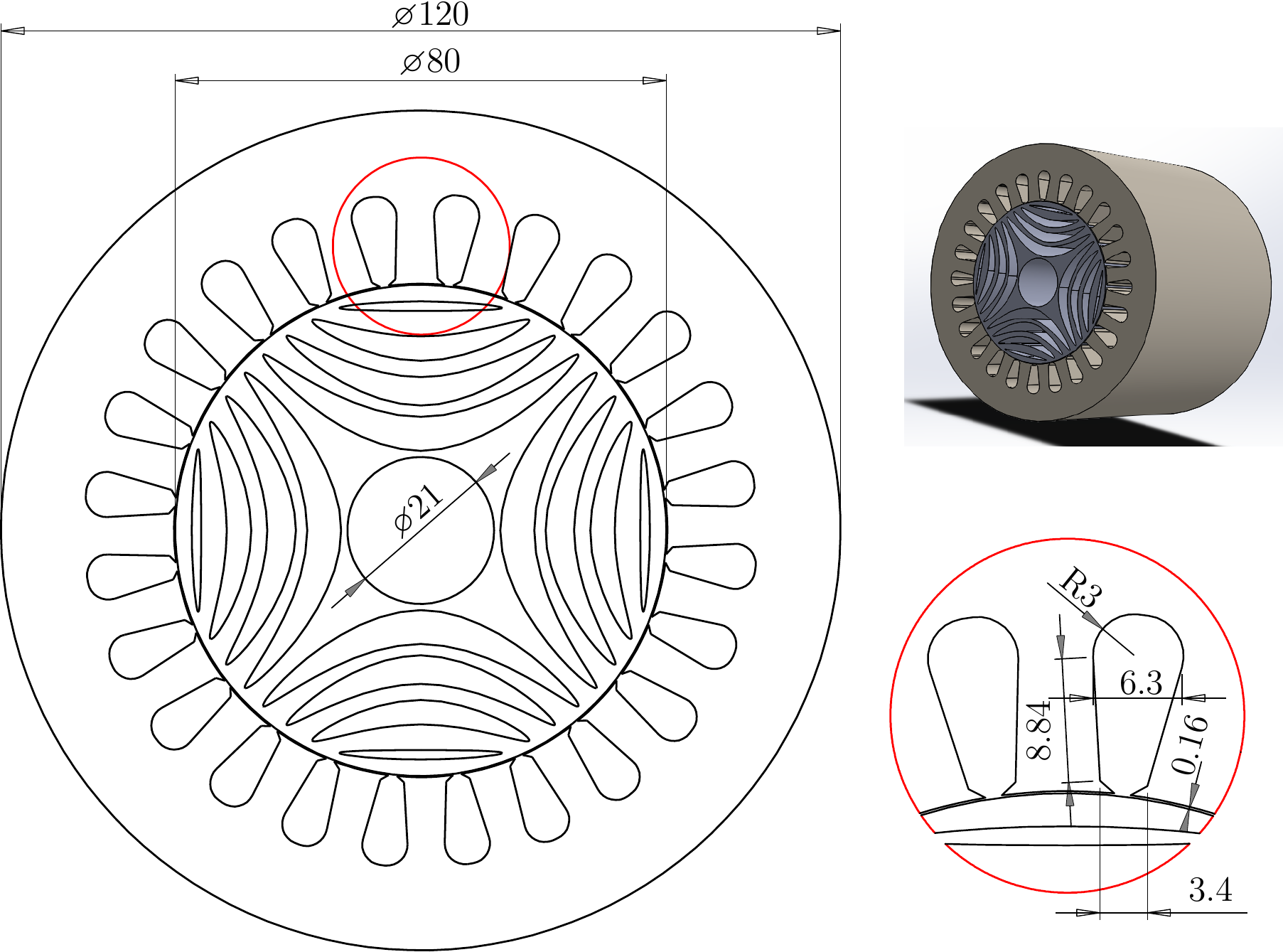}
	\caption{Reluctance motor geometry}
	\label{fig:motor_obr}
\end{figure}

The goal of this project was to perform a shape optimization on a reluctance machine (see Fig. \ref{fig:motor_obr}). The shape of the machine was described by Bézier curves. The position of the control points were optimized. 
The physics of the problem was described by the partial differential equations of the magnetic field:
\begin{eqnarray}
  \gamma \frac{\partial \vec{A}}{\partial t} + \curl \vec{H} - \gamma \vec{v} \times \vec{B} &=& \vec{J}_{\rm e}, \\
\curl \vec{A} = \vec{B},\quad \vec{B}=\mu \vec{H}\,,  
\end{eqnarray}
where $\vec{A}$ is the magnetic vector potential, $\vec{H}$ is the magnetic field strength, $\vec{B}$ is the magnetic flux density, $\vec{v}$ is velocity, symbol $\omega$ means for frequency, $\gamma$ stands for conductivity, $\mu$ for magnetic permeability and $J_\mathrm{e}$ is the current density.

\begin{figure}[!ht]
	\subfloat[]{\includegraphics[width=.48 \linewidth]{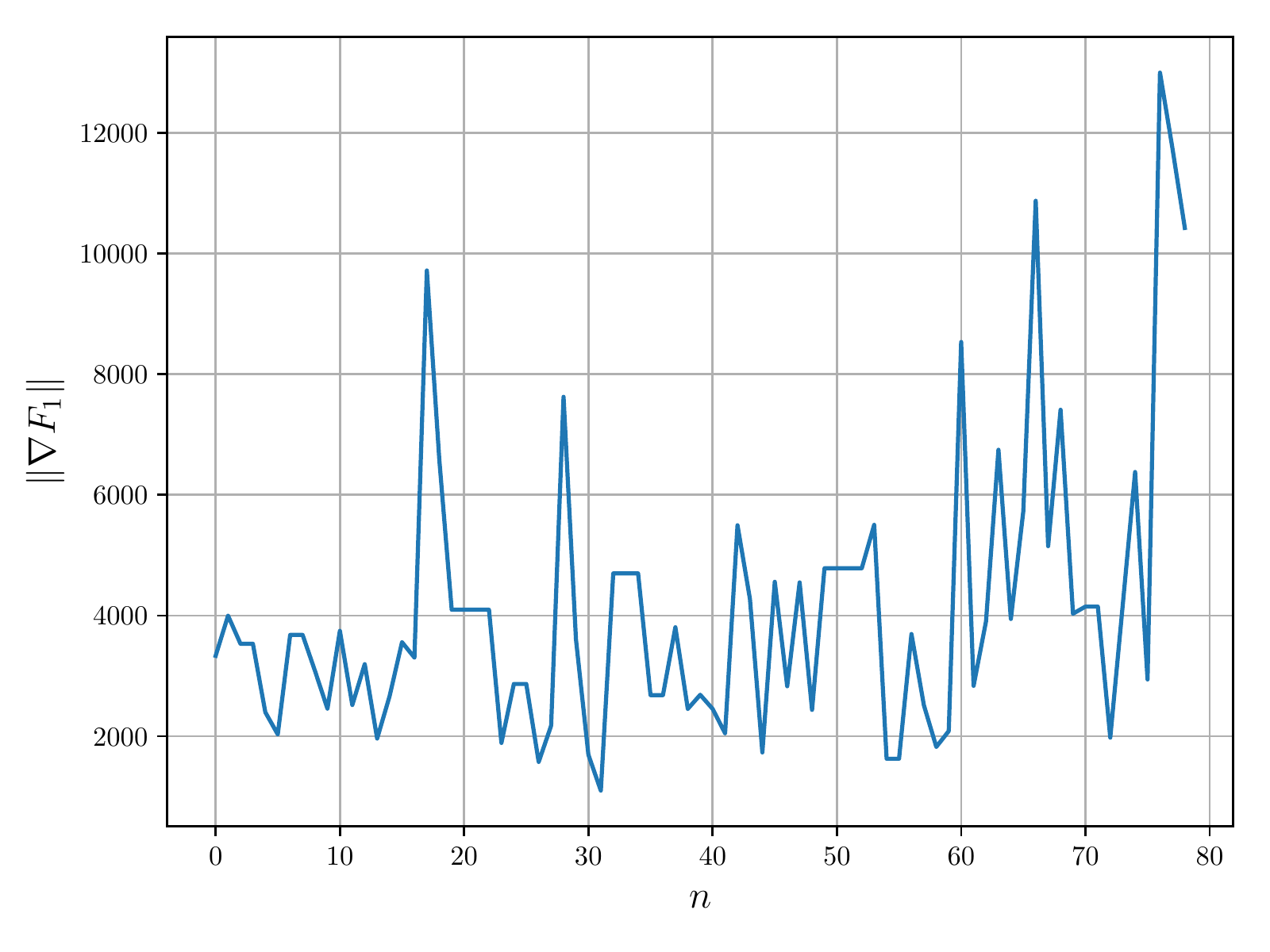} \label{fig:grad_F1}}
	\subfloat[]{\includegraphics[width=.48 \linewidth]{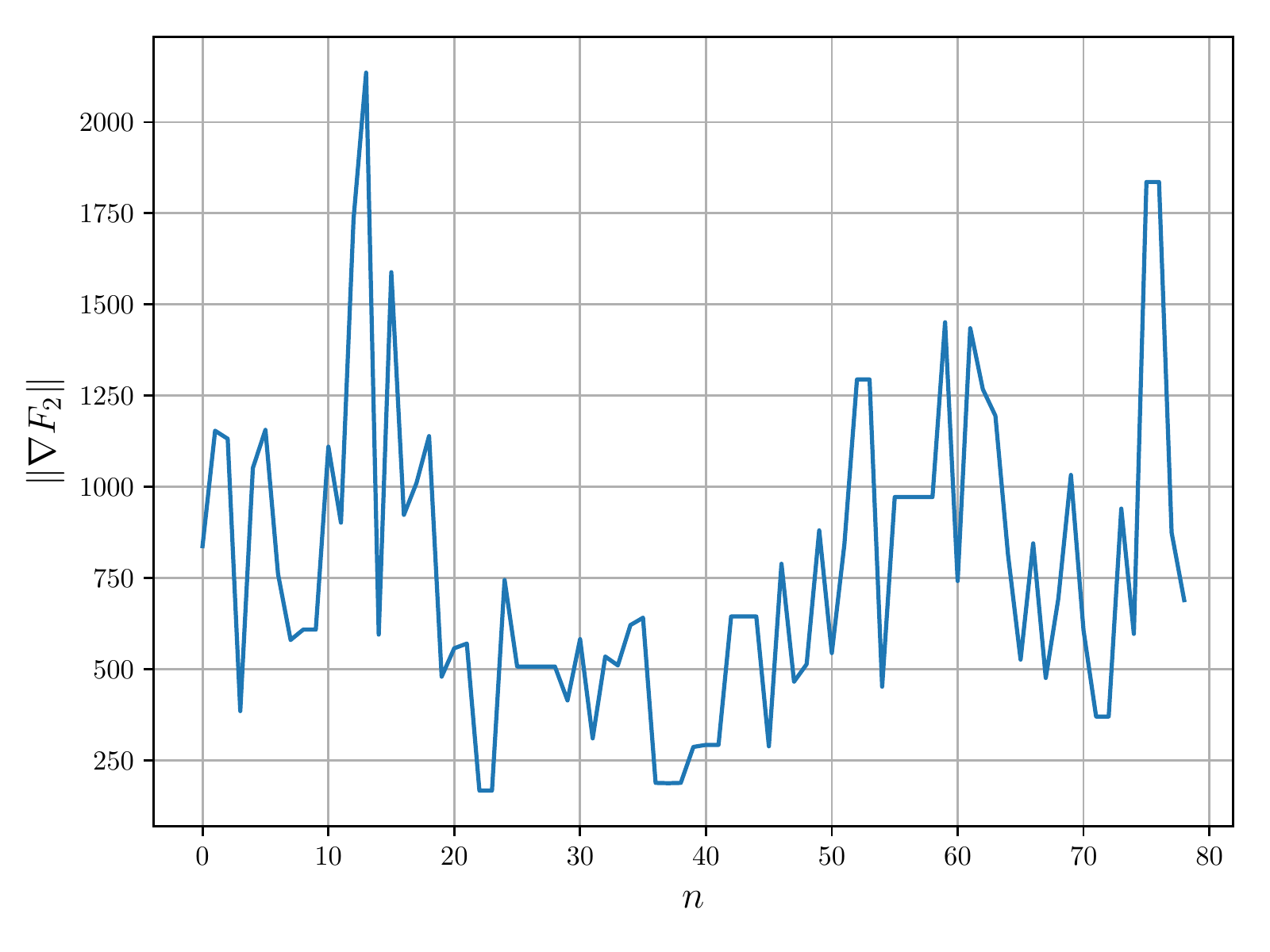} \label{fig:grad_F2}} 
	\caption{Norm of the gradient of the goal function $F_1$-(a), and the goal function $F_2$-(b).}
    \label{fig:norm}
\end{figure}

There were two objective functions considered. The first objective ($F_1$) was the mechanical torque: 
\begin{equation}
\vec{T}_{ \rm M} = \int_S \vec{r}\times\left(\vec{S}_{\rm M}\cdot\dif {\vec S}\right)\,,
\label{eq:maxwell_stress_tensor_torque}
\end{equation}
where $\vec{r}$ is the position vector and 
\begin{equation}
\vec{S}_{\rm M} = -\frac{1}{2\mu} (\vec{B} \cdot \vec{B}) \mat{I} + \frac{1}{\mu} \vec{B} \otimes \vec{B}.
\end{equation}
is the Maxwell stress tensor.
The second objective ($F_2$) was the standard deviation of the mechanical torque. NSGA-II algorithm was selected for the optimization. The evolution of the resulted populations is depicted in Fig. \ref{fig:motor_obr}, the individuals of the Pareto-front are highlighted by black color. It can be seen from the results, that each individual have different sensitivity on the parameters. Because of the manufacturing process works with limited precision, to ensure a robust design, it is necessary to select those individuals which has the lowest sensitivity. As a measure of the sensitivity, the norm of the goal function's gradient was selected (Fig. \ref{fig:norm}). 
\begin{equation}
    s_1 = \| \grad(F_1) \|,  s_1 = \| \grad(F_2) \|.
\end{equation}

\begin{figure}[!ht]
	\centering
	\includegraphics[width=.99 \linewidth]{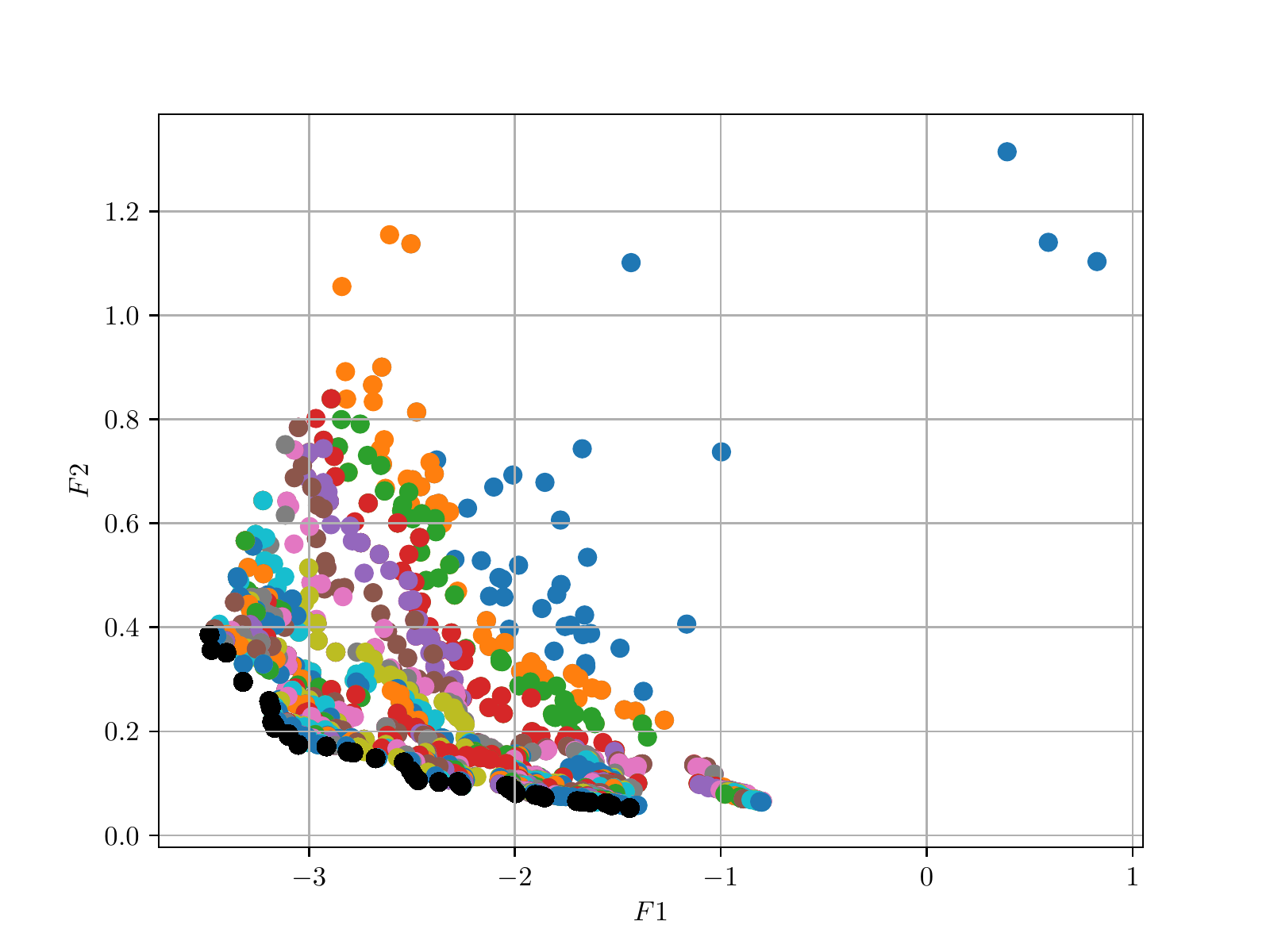}
	\caption{Evolution of goal functions during the shape optimization.}
	\label{fig:motor_optimization}
\end{figure}

\subsection{TEAM Benchmark}

Nowadays, much effort is devoted to accelerate the numerical processes in optimization and inverse techniques in industrial and other domains. This example is aimed to solve a seemingly simple, benchmark problem, to find out an optimal distribution of current-carrying turns that should generate a homogeneous magnetic field in a given region \cite{di2018benchmark}.

\begin{figure}[!ht]
\centering
\includegraphics[width=0.5\linewidth]{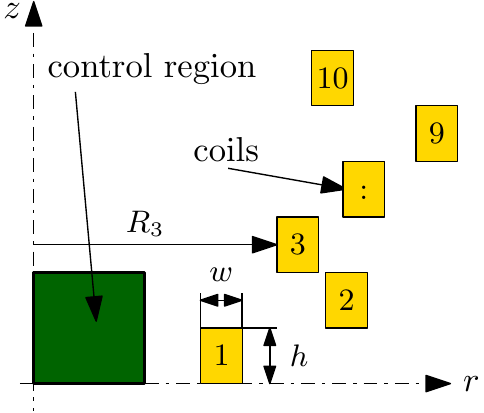} 
\caption{Geometry and design variables}
\label{fig:geometry_benchmark}
\end{figure}

The arrangement of the system is depicted in Fig.\,\ref{fig:geometry_benchmark}    

The geometry is known except for the radii of particular turns that may vary within the prescribed ranges. The values of direct currents in the turns and dimensions of the control region are also given. The task is to find such radii of particular turns that would produce as uniform magnetic field in the indicated region as possible.

Accurate semi-analytical solution can be derived. The components $B_r$ and $B_z$ of magnetic flux density produced at a point $P(R,Z)$ are given by the formulae
$$
B_{r}(R,Z)=C\cdot[g(R_2,R,Z_2-Z)-g(R_2,R,Z_1-Z)-
$$
\begin{equation}
-g(R_1,R,Z_2-Z)+g(R_1,R,Z_1-Z)]
\end{equation}
and
$$
B_{z}(R,Z)=C\cdot[h(R_2,R,Z_2-Z)-h(R_2,R,Z_1-Z)-
$$
\begin{equation}
-h(R_1,R,Z_2-Z)+h(R_1,R,Z_1-Z)]\,,
\end{equation}
where
\begin{equation}
C=\frac{\mu_0 I}{4\pi (Z_2-Z_1)\ln(R_2/R_1)}\,.
\end{equation}

Here, for example
\begin{equation}
g(R_2,R,Z_2-Z)=\int_{\varphi=0}^{2\pi}\ln[R_2-R\cdot\cos\varphi+d_{22}]\cos\varphi\,{\rm d}\varphi\,,
\end{equation}
\begin{equation}
h(R_2,R,Z_2-Z)=-\int_{\varphi=0}^{2\pi}\ln[Z_2-Z+d_{22}]\,{\rm d}\varphi\,,
\end{equation}
and
\begin{equation}
d_{22}=\sqrt{R_2^2+R^2-2R_2 R\cos\varphi+(Z-Z_2)^2}\,.
\end{equation}

The other functions are obtained by standard interchanging of the indices. The last integral with respect to $\varphi$ is calculated using the Gauss quadrature formulae. Magnetic field produced by more turns is then given by the superposition of the partial fields.

The inverse task is solved by Artap. The main ideas of the Gaussian Process Regressor are described in \cite{gaussian}. After $N$ standard iterations (in our case $N=30$) a surrogate model is created based on the Gaussian Process Regressor. Prior to further evaluations of the objective function (which holds both for single-criterion and multi-criteria processes), an estimate is carried out, which provides the estimated value and standard deviation. If this deviation is low, the above estimated value is taken as the new value of the objective function, otherwise this function is evaluated in the standard way. This step is then repeated after further $N$ iteration, which makes the surrogate model more precise. The optimization process is now substantially faster. 

The goal function for single-objective problem is to design the geometry of coils that minimize discrepancy between the prescribed valued $\vec{B}_0$ and actual distribution of $\vec{B}$ in the region of interest. 
$$
F_1(r) = \sup_{q=1,n_{\rm p}} |\vec{B}(r_q,z_q)-\vec{B}_0(r_q,z_q)|,
$$
where $B_0(r_q,z_q) = 0.2\,{\rm mT}$ is the prescribed value of magnetic flux density and $n_{\rm p}$ is the number of points.
In the multi-objective problem, the following sensitivity function is minimized
$$
F_2(r) = \sup_{q=1,n_{\rm p}} (|\vec{B}(r_q,z_q)-\vec{B}(r_q+\Delta r,z_q)| + 
$$
$$
+ |\vec{B}(r_q,z_q)-\vec{B}(r_q-\Delta r,z_q)|).
$$

For the multi-objective optimization, the NSGA-II algorithm was used. The number of populations is 80 and population size is 20.
\begin{table}[!ht]
\renewcommand{\arraystretch}{1.3}
\caption{Multi-objective optimization (NSGA-II)}
\label{tab:single}
\centering
\begin{tabular}{|l||c|c|c|}
\hline
  & Elapsed time & Predict. & Eval. \\
\hline
Surrogate & 341 s & 1243 & 327 \\
\hline
Evaluation & 1091 s & 0 & 1091 \\
\hline
\end{tabular}
\end{table}

\begin{figure}[!ht]
\centering
\includegraphics[width=0.9\linewidth]{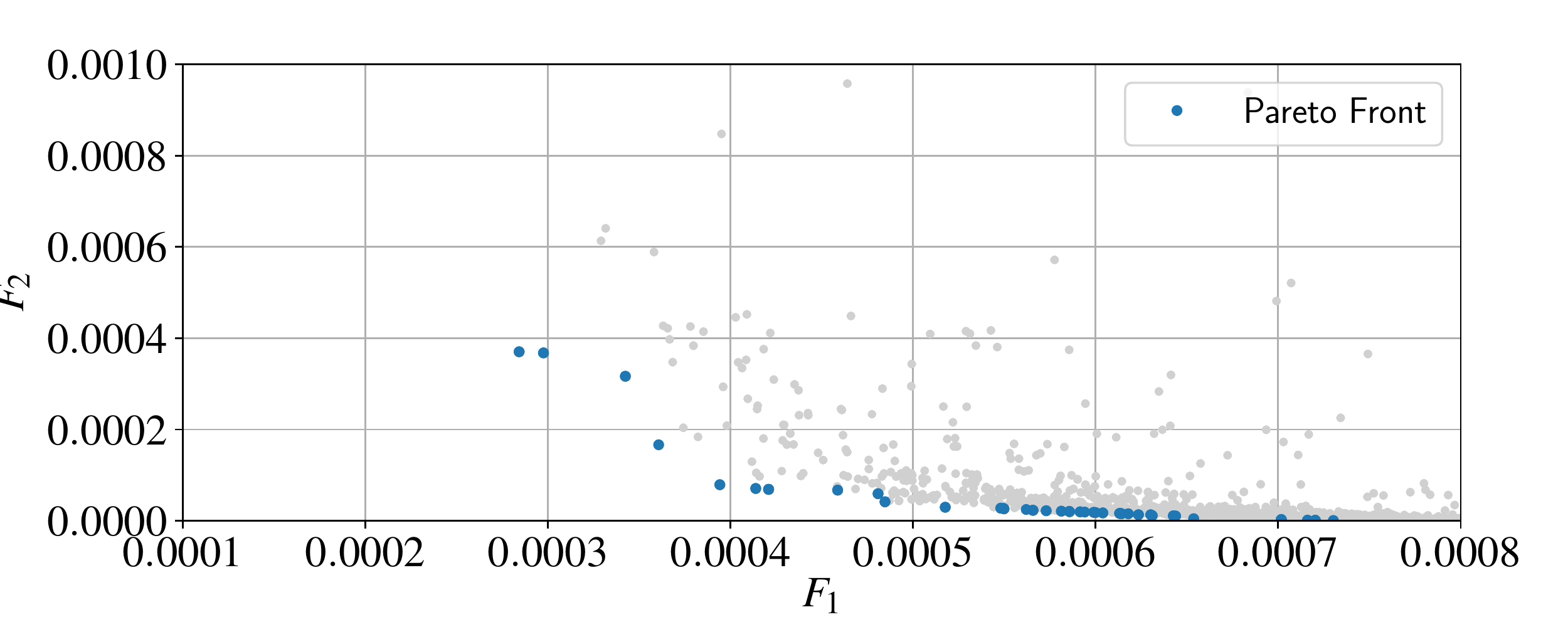}
\caption{Pareto front - without acceleration}
\label{fig:nsgaii_f1f3}
\end{figure}

\begin{figure}[!ht]
\centering
\includegraphics[width=0.9\linewidth]{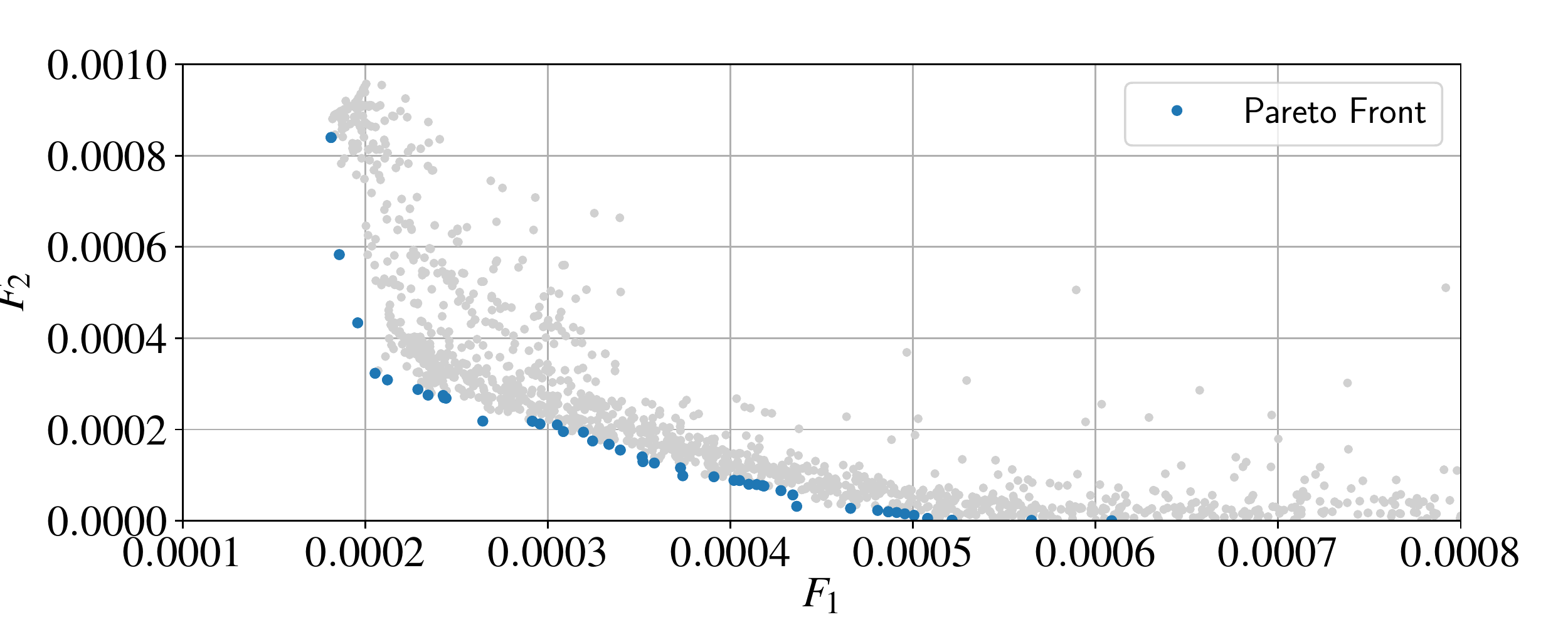}
\caption{Pareto front - with surrogate regression}
\label{fig:nsgaii_f1f3}
\end{figure}

\section{Conclusion}

In this paper, Artap, an open source high-performance computing platform for multidisciplinary, especially robust design optimization, has been presented. 
The applicability of the code  has been demonstrated on an induction brazing application. 
Artap facilitates the use of many different numerical tools together during the calculation process. 
This is essential to solve a robust design optimization problem and ensures a straightforward way to solve other complex inverse problems, like parameter identification, model order reduction and shape optimization problems, as well. 
Due to the three-layered architecture of Artap, the solvers can be changed flexibly and allows to try the most suitable solver for the design optimization problem. 
Future work on Artap -- besides the extension with additional optimization algorithms and numerical solvers -- will give a simple, web-based interface to the project.

Artap is a MIT licensed Software, which is available for download from the Agros-suite website:\\
http://www.agros2d.org/artap/

\section*{Acknowledgement}
This research has been supported by the Ministry of Education, Youth and Sports of the Czech Republic under the RICE New Technologies and Concepts for Smart Industrial Systems, project No. LO1607 and by an internal project SGS-2018-043.

\bibliographystyle{IEEEtran}
\bibliography{references}

\begin{thebibliography}{10}
\providecommand{\url}[1]{#1}
\csname url@samestyle\endcsname
\providecommand{\newblock}{\relax}
\providecommand{\bibinfo}[2]{#2}
\providecommand{\BIBentrySTDinterwordspacing}{\spaceskip=0pt\relax}
\providecommand{\BIBentryALTinterwordstretchfactor}{4}
\providecommand{\BIBentryALTinterwordspacing}{\spaceskip=\fontdimen2\font plus
\BIBentryALTinterwordstretchfactor\fontdimen3\font minus
  \fontdimen4\font\relax}
\providecommand{\BIBforeignlanguage}[2]{{%
\expandafter\ifx\csname l@#1\endcsname\relax
\typeout{** WARNING: IEEEtran.bst: No hyphenation pattern has been}%
\typeout{** loaded for the language `#1'. Using the pattern for}%
\typeout{** the default language instead.}%
\else
\language=\csname l@#1\endcsname
\fi
#2}}
\providecommand{\BIBdecl}{\relax}
\BIBdecl

\bibitem{Ben-Tal2002}
\BIBentryALTinterwordspacing
A.~Ben-Tal and A.~Nemirovski, ``Robust optimization -- methodology and
  applications,'' \emph{Mathematical Programming}, vol.~92, no.~3, pp.
  453--480, May 2002. [Online]. Available:
  \url{https://doi.org/10.1007/s101070100286}
\BIBentrySTDinterwordspacing

\bibitem{li2005multi}
M.~Li, S.~Azarm, and V.~Aute, ``A multi-objective genetic algorithm for robust
  design optimization,'' in \emph{Proceedings of the 7th annual conference on
  Genetic and evolutionary computation}.\hskip 1em plus 0.5em minus 0.4em\relax
  ACM, 2005, pp. 771--778.

\bibitem{tsui1992overview}
K.-L. Tsui, ``An overview of taguchi method and newly developed statistical
  methods for robust design,'' \emph{Iie Transactions}, vol.~24, no.~5, pp.
  44--57, 1992.

\bibitem{rdo}
\BIBentryALTinterwordspacing
I.~Doltsinis and Z.~Kang, ``Robust design of structures using optimization
  methods,'' \emph{Computer Methods in Applied Mechanics and Engineering}, vol.
  193, no.~23, pp. 2221 -- 2237, 2004. [Online]. Available:
  \url{http://www.sciencedirect.com/science/article/pii/S0045782504000787}
\BIBentrySTDinterwordspacing

\bibitem{ih_review}
O.~{Luc\'{i}a}, P.~{Maussion}, E.~J. {Dede}, and J.~M. {Burd\'{i}o},
  ``Induction heating technology and its applications: Past developments,
  current technology, and future challenges,'' \emph{IEEE Transactions on
  Industrial Electronics}, vol.~61, no.~5, pp. 2509--2520, May 2014.

\bibitem{induction_brazing_1}
D.~P\'anek, P.~Karban, and I.~Dole\v{z}el, ``Comparison of simplified
  techniques for solving selected coupled electroheat problems,'' to appear in
  Proceedings of the HES-19, International Symposium on Heating by
  Electromagnetic Sources, 2019, p.~2p.

\bibitem{induction_brazing_2}
D.~P\'anek, T.~Orosz, P.~Krop\'ik, P.~Karban, and I.~Dole\v{z}el,
  ``Reduced-order model based temperature control of induction brazing
  process,'' to appear in 2019 Electric Power Quality and Supply Reliability
  (PQ), 2019, p.~4p.

\bibitem{8124818}
K.~{Daukaev}, A.~{Rassõlkin}, A.~{Kallaste}, T.~{Vaimann}, and A.~{Belahcen},
  ``A review of electrical machine design processes from the standpoint of
  software selection,'' in \emph{2017 IEEE 58th International Scientific
  Conference on Power and Electrical Engineering of Riga Technical University
  (RTUCON)}, Oct 2017, pp. 1--6.

\bibitem{6479303}
Y.~{Duan} and D.~M. {Ionel}, ``A review of recent developments in electrical
  machine design optimization methods with a permanent-magnet synchronous motor
  benchmark study,'' \emph{IEEE Transactions on Industry Applications},
  vol.~49, no.~3, pp. 1268--1275, May 2013.

\bibitem{metaheuristic}
T.~Orosz, A.~Sleisz, and Z.~A. Tamus, ``Metaheuristic optimization preliminary
  design process of core-form autotransformers,'' \emph{IEEE Transactions on
  Magnetics}, vol.~52, no.~4, pp. 1--10, April 2016.

\bibitem{trafoevo}
T.~Orosz, ``Evolution and modern approaches of the power transformer cost
  optimization methods,'' \emph{Periodica Polytechnica Electrical Engineering
  and Computer Science}, 2019, to be published.

\bibitem{PAVLICEK2019}
\BIBentryALTinterwordspacing
K.~Pavlíček, V.~Kotlan, and I.~Doležel, ``Applicability and comparison of
  surrogate techniques for modeling of selected heating problems,''
  \emph{Computers \& Mathematics with Applications}, 2019. [Online]. Available:
  \url{http://www.sciencedirect.com/science/article/pii/S0898122119300811}
\BIBentrySTDinterwordspacing

\bibitem{moulard2014software}
T.~Moulard, B.~Chr{\'e}tien, and E.~Yoshida, ``Software tools for nonlinear
  optimization,'' \emph{J. Robot. Soc. Jpn.}, vol.~32, no.~6, pp. 536--541,
  2014.

\bibitem{moulard2013roboptim}
T.~Moulard, F.~Lamiraux, K.~Bouyarmane, and E.~Yoshida, ``Roboptim: an
  optimization framework for robotics,'' in \emph{Robomec}, 2013, p.~4p.

\bibitem{biscani2010global}
F.~Biscani, D.~Izzo, and C.~H. Yam, ``A global optimisation toolbox for
  massively parallel engineering optimisation,'' \emph{arXiv preprint
  arXiv:1004.3824}, 2010.

\bibitem{Andersson2019}
\BIBentryALTinterwordspacing
J.~A.~E. Andersson, J.~Gillis, G.~Horn, J.~B. Rawlings, and M.~Diehl, ``Casadi:
  a software framework for nonlinear optimization and optimal control,''
  \emph{Mathematical Programming Computation}, vol.~11, no.~1, pp. 1--36, Mar
  2019. [Online]. Available: \url{https://doi.org/10.1007/s12532-018-0139-4}
\BIBentrySTDinterwordspacing

\bibitem{noauthor_openmdao.org_nodate}
\BIBentryALTinterwordspacing
``\BIBforeignlanguage{en-US}{{OpenMDAO}.org {\textbar} {An} open-source
  framework for efficient multidisciplinary optimization.{OpenMDAO}.org
  {\textbar} {An} open-source framework for efficient multidisciplinary
  optimization.}'' [Online]. Available: \url{http://openmdao.org/}
\BIBentrySTDinterwordspacing

\bibitem{noauthor_platypus_nodate}
\BIBentryALTinterwordspacing
``Platypus - {Multiobjective} {Optimization} in {Python} — {Platypus}
  documentation.'' [Online]. Available:
  \url{https://platypus.readthedocs.io/en/latest/}
\BIBentrySTDinterwordspacing

\bibitem{johnson2014nlopt}
S.~G. Johnson, ``The nlopt nonlinear-optimization package,'' 2014.

\bibitem{jones2014scipy}
E.~Jones, T.~Oliphant, and P.~Peterson, ``$\{$SciPy$\}$: Open source scientific
  tools for $\{$Python$\}$,'' 2014.

\bibitem{toth2018dual}
B.~T{\'o}th, ``Dual and mixed nonsymmetric stress-based variational
  formulations for coupled thermoelastodynamics with second sound effect,''
  \emph{Continuum Mechanics and Thermodynamics}, vol.~30, no.~2, pp. 319--345,
  2018.

\bibitem{bangerth2007deal}
W.~Bangerth, R.~Hartmann, and G.~Kanschat, ``deal. ii—a general-purpose
  object-oriented finite element library,'' \emph{ACM Transactions on
  Mathematical Software (TOMS)}, vol.~33, no.~4, p.~24, 2007.

\bibitem{karban_numerical_2013}
\BIBentryALTinterwordspacing
P.~Karban, F.~Mach, P.~Kůs, D.~Pánek, and I.~Doležel,
  ``\BIBforeignlanguage{en}{Numerical solution of coupled problems using code
  {Agros}2d},'' \emph{\BIBforeignlanguage{en}{Computing}}, vol.~95, no.~1, pp.
  381--408, May 2013. [Online]. Available:
  \url{https://link.springer.com/article/10.1007/s00607-013-0294-4}
\BIBentrySTDinterwordspacing

\bibitem{deb2002fast}
K.~Deb, A.~Pratap, S.~Agarwal, and T.~Meyarivan, ``A fast and elitist
  multiobjective genetic algorithm: Nsga-ii,'' \emph{IEEE transactions on
  evolutionary computation}, vol.~6, no.~2, pp. 182--197, 2002.

\bibitem{coello2007evolutionary}
C.~A.~C. Coello, G.~B. Lamont, D.~A. Van~Veldhuizen \emph{et~al.},
  \emph{Evolutionary algorithms for solving multi-objective problems}.\hskip
  1em plus 0.5em minus 0.4em\relax Springer, 2007, vol.~5.

\bibitem{mishra2002fast}
S.~K. Mishra, P.~Ganapati, S.~Meher, and R.~Majhi, ``A fast multiobjective
  evolutionary algorithm for finding wellspread pareto-optimal solutions,'' in
  \emph{In KanGAL Report No. 2003002, lndian Institute Of Technology
  Kanpur}.\hskip 1em plus 0.5em minus 0.4em\relax Citeseer, 2002.

\bibitem{kennedy1995eberhart}
J.~Kennedy, ``Eberhart, r.: Particle swarm optimization,'' in \emph{Proceedings
  of IEEE international conference on neural networks}, vol.~4, no.~2.\hskip
  1em plus 0.5em minus 0.4em\relax IEEE Press, 1995, pp. 1942--1948.

\bibitem{gao2012implementing}
F.~Gao and L.~Han, ``Implementing the nelder-mead simplex algorithm with
  adaptive parameters,'' \emph{Computational Optimization and Applications},
  vol.~51, no.~1, pp. 259--277, 2012.

\bibitem{powell1964efficient}
M.~J. Powell, ``An efficient method for finding the minimum of a function of
  several variables without calculating derivatives,'' \emph{The computer
  journal}, vol.~7, no.~2, pp. 155--162, 1964.

\bibitem{nocedal2006sequential}
J.~Nocedal and S.~J. Wright, ``Sequential quadratic programming,''
  \emph{Numerical optimization}, pp. 529--562, 2006.

\bibitem{lenders2018trlib}
F.~Lenders, C.~Kirches, and A.~Potschka, ``trlib: A vector-free implementation
  of the gltr method for iterative solution of the trust region problem,''
  \emph{Optimization Methods and Software}, vol.~33, no.~3, pp. 420--449, 2018.

\bibitem{conn2000trust}
A.~R. Conn, N.~I. Gould, and P.~L. Toint, \emph{Trust region methods}.\hskip
  1em plus 0.5em minus 0.4em\relax Siam, 2000, vol.~1.

\bibitem{martinez2014bayesopt}
R.~Martinez-Cantin, ``Bayesopt: A bayesian optimization library for nonlinear
  optimization, experimental design and bandits,'' \emph{The Journal of Machine
  Learning Research}, vol.~15, no.~1, pp. 3735--3739, 2014.

\bibitem{pedregosa2011scikit}
F.~Pedregosa, G.~Varoquaux, A.~Gramfort, V.~Michel, B.~Thirion, O.~Grisel,
  M.~Blondel, P.~Prettenhofer, R.~Weiss, V.~Dubourg \emph{et~al.},
  ``Scikit-learn: Machine learning in python,'' \emph{Journal of machine
  learning research}, vol.~12, no. Oct, pp. 2825--2830, 2011.

\bibitem{di2018benchmark}
P.~Di~Barba, M.~E. Mognaschi, D.~A. Lowther, and J.~K. Sykulski, ``A benchmark
  team problem for multi-objective pareto optimization of electromagnetic
  devices,'' \emph{IEEE Transactions on Magnetics}, vol.~54, no.~3, pp. 1--4,
  2018.

\bibitem{chatterjee2000introduction}
A.~Chatterjee, ``An introduction to the proper orthogonal decomposition,''
  \emph{Current science}, pp. 808--817, 2000.

\bibitem{gaussian}
C.~E. Rasmussen and C.~K.~I. Williams, \emph{Gaussian Processes for Machine
  Learning}.\hskip 1em plus 0.5em minus 0.4em\relax MIT Press, 2006.

\end{thebibliography}

\end{document}